\documentclass[12pt]{article}
\usepackage{geometry}
\usepackage[utf8]{inputenc}
\usepackage{amsmath, amsfonts, amsthm, amssymb}
\usepackage{enumerate}
\usepackage{mathtools}
\usepackage[hidelinks]{hyperref}
\usepackage[noabbrev]{cleveref}
\usepackage{tabu}
\usepackage{algorithm}
\usepackage{algorithmic}
\usepackage{mathrsfs}
\usepackage{xcolor}
\usepackage{caption}
\usepackage{subcaption}
\usepackage{tikz-cd}
\usepackage{xspace}
\usepackage{amsrefs}

\usepackage{colortbl}

\newcommand{\im}{\textup{Im}}
\newcommand{\st}{\,:\,}
\renewcommand{\mod}{\ \textup{mod}\,}

\newcommand{\F}{\mathbb{F}}

\newcommand{\zm}{\mathbb{Z}^m}
\newcommand{\z}{\mathbb{Z}}

\newcommand{\aut}{\textup{Aut}}

\allowdisplaybreaks

\DeclarePairedDelimiter\set{\{}{\}}
\DeclarePairedDelimiter\abs{\lvert}{\rvert}
\DeclarePairedDelimiter\norm{\lVert}{\rVert}
\DeclarePairedDelimiter\inprod{\langle}{\rangle}

\makeatletter 
\let\oldset\set
\def\set{\@ifstar{\oldset}{\oldset*}}

\let\oldabs\abs
\def\abs{\@ifstar{\oldabs}{\oldabs*}}

\let\oldnorm\norm
\def\norm{\@ifstar{\oldnorm}{\oldnorm*}}

\let\oldinprod\inprod
\def\inprod{\@ifstar{\oldinprod}{\oldinprod*}}
\makeatother

\newtheorem{theorem}{Theorem}
\newtheorem{lemma}[theorem]{Lemma}
\newtheorem{conjecture}[theorem]{Conjecture}
\newtheorem{cor}[theorem]{Corollary}
\newtheorem{prop}[theorem]{Proposition}

\theoremstyle{definition}
\newtheorem{definition}[theorem]{Definition}

\newtheorem{example}[theorem]{Example}

\newtheorem{remark}[theorem]{Remark}

\title{Circular Costas maps: a multidimensional analog of circular Costas sequences}
\author{Ivelisse Rubio and Jaziel Torres}

\begin{document}

\maketitle

\begin{abstract}
    A unifying theoretical framework is presented, in which the connections among Costas sequences, circular Costas sequences, Costas polynomials, the shifting property, and Welch sequences are extended to the multidimensional context.    Several conjectures on multidimensional periodic Costas arrays by J. Ortiz-Ubarri et al. are proved.
    Furthermore, a conjecture on Costas polynomials over extension fields presented by Muratovi\'c-Ribi\'c et al. is showed to be a multidimensional extension of a conjecture by Golomb and Moreno on circular Costas sequences.
    A weaker version of said conjecture is proved by considering a multidimensional extension of the shifting Costas property defined by O. Moreno.
\end{abstract}

\section{Introduction}
In the quest to improve the performance of SONAR detection, in a 1965 technical report \cite{costas1965medium}, later published as a journal article \cite{costas1984study}, John P. Costas proposed frequency hop patterns represented as arrays of zeros and ones that can unambiguously determine the distance and velocity of a moving target.
These arrays are constructed as follows:
\begin{quote}
    Place $n$ ones in an otherwise null $n\times n$ matrix such that each row contains a single one as does each column. 
    Make the placement such that for all possible $x$-$y$ shift combinations of the resulting (permutation) matrix relative to itself, at most one pair of ones will coincide. \cite{costas1984study}*{p.~997}
\end{quote}

The arrays proposed by J. P. Costas were originally called \emph{constellation arrays} \citelist{\cite{cohen1991primitive}*{Section I} \cite{golomb1982two}*{Section II}}, but later acquired the name for which they are known today: \emph{Costas arrays}.

Using the usual notation $[n] = \set{1,2,\dots,n}$, a permutation matrix of order $n$ is equivalent to a bijection $\varphi:[n] \to [n]$, where $\varphi(i) = t$ if and only if there is a one in position $(i,t)$.
Given such bijection, one can consider the sequence of images $\varphi(1), \varphi(2), \dots, \varphi(n)$.
Thus, Costas arrays can be equivalently described as sequences.

\begin{definition}\label{def:costas_seq}
	A sequence $a_0, a_1, \dots, a_{n-1}$ which is a permutation of the integers $1, 2, \dots, n$ is a \textbf{Costas sequence} if it satisfies: for any $k\in\set{1, 2, \dots, n-1}$,
	\begin{equation}\label{eq:costas_distinct_diff}
		a_{i+k}-a_i = a_{j+k}-a_j \implies i=j 
	\end{equation}
	for all $i,j$ with $0\leqslant i, j \leqslant n-1-k$.
\end{definition}
Condition \eqref{eq:costas_distinct_diff} is called the \emph{distinct difference property}.
A permutation sequence has the distinct difference property if and only if it defines a Costas array.

\begin{figure}[ht]
    \centering
    \includegraphics[scale=.2]{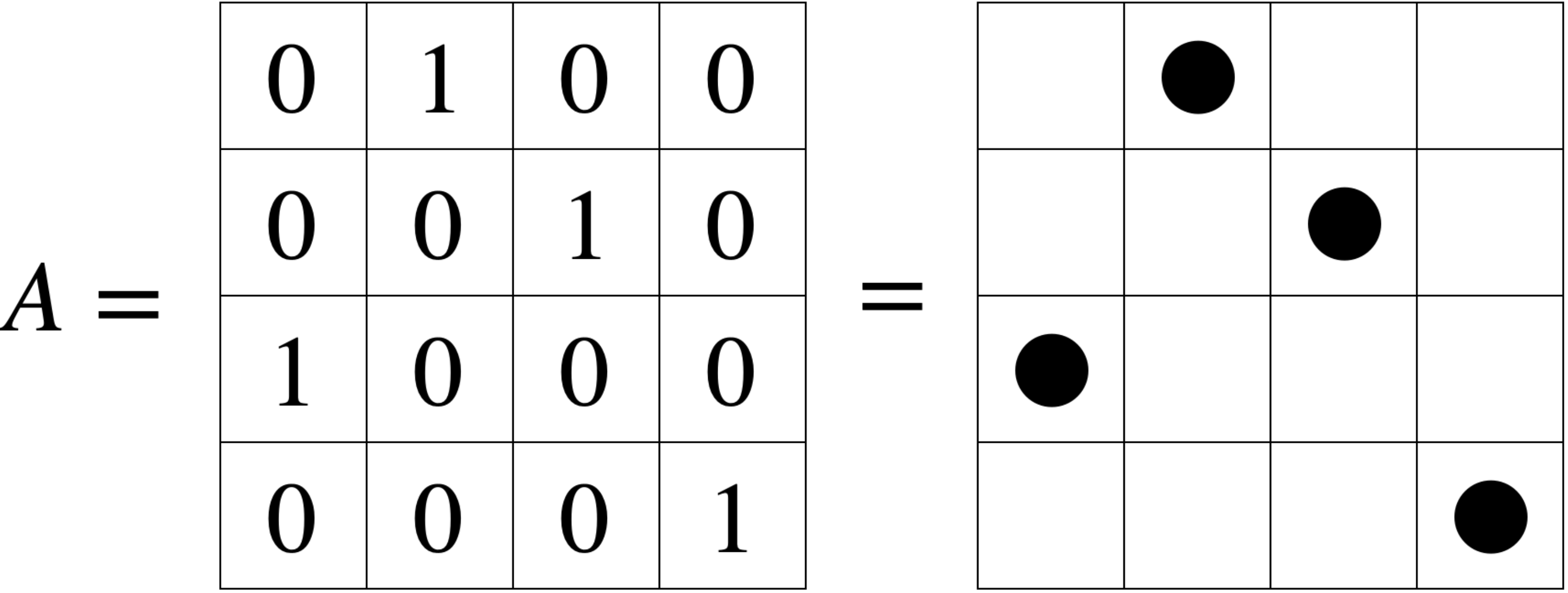}
    \caption{A Costas array of order $4$ corresponding to the sequence $2,4,3,1$.}
\end{figure}

One can verify the distinct difference property by computing the sequence's difference triangle, where the $k$-th row of the triangle contains $a_{i+k}-a_i$ for $i = 0, 1, \dots, n-1-k$, where $1 \leqslant k \leqslant n-1$, as exemplified in \Cref{table}.
The sequence is Costas if and only if each row has no repeats.

\begin{table}
    \centering
    \begin{tabular}{l||rrrr}
    $k$ & \textbf{2} & \textbf{4} & \textbf{3} & \textbf{1}\\ \hline
    \textit{1} & 2 & $-1$ & $-2$ & \hspace{10pt} \ \\
    \textit{2} & 1 & $-3$ & & \\
    \textit{3} & $-1$ & & & 
    \end{tabular}
    \caption{Difference triangle for the Costas sequence $2,4,3,1$.}
    \label{table}
\end{table}

One of the few systematic procedures to construct Costas arrays is attributed to Lloyd R. Welch \cite{golomb1982two}, and is as follows:

\begin{prop}[Welch construction]\label{prop:welch}
    Let $\alpha$ be a primitive element in $\F_p$, with $p$ prime, and $c\in\z$.
    The sequence $\alpha^{0+c}, \alpha^{1 + c}, \dots, \alpha^{p-2 + c}$ is a Costas sequence.  
\end{prop}

In a 1984 article, Solomon W. Golomb and Herbert Taylor considered some periodic properties of Costas arrays \cite{golomb-taylor1984}*{Section III.A}.
In particular, they mentioned:
\begin{quote}
	Repeating the $2\times 2$ Costas array in both directions over the entire plane gives a doubly periodic checkerboard pattern with a Costas array in every $2\times 2$ window. 
	For any $n > 2$, however, there does not exist a doubly periodic pattern with a Costas array in every $n \times n$ window.
	\cite{golomb-taylor1984}*{p.~1154}
\end{quote}
The above result was proved by Taylor in \cite{taylor1984non}*{Theorem B}.

On the other hand, singly periodical patterns do exist.
That is, repeating an $n \times n$ array horizontally to obtain an array with $n$ rows and infinite columns such that every $n \times n$ window is a Costas array.
In the context of single periodicity, Taylor's theorem \cite{taylor1984non}*{Theorem B} states that, if an $n \times n$ Costas array is singly periodical and its transpose is also singly periodical, then $n = 2$;
an array too small to be interesting, let alone useful.
Equivalently, a Costas sequence is singly periodical if any circular shift of its terms is also a Costas sequence.
However, the only known singly periodical sequences are those obtained from the Welch construction, which lead Golomb and Taylor to conjecture that those are indeed the only ones (problem \#8 in \cite{golomb-taylor1984}*{Section V}).
This conjecture, published in 1984, is still open:
\begin{conjecture}[Golomb \& Taylor \cite{golomb-taylor1984}]\label{conj:singly_periodic}
	A Costas sequence is singly periodical if and only if it is Welch.
\end{conjecture}

Making progress towards \Cref{conj:singly_periodic}, in 1992, Oscar Moreno \cite{moreno1992shifting} proved that a condition he called the \emph{shifting property}, which is stronger than single periodicity, characterizes the Welch construction.

\begin{definition}\label{def:moreno_shifting}
    For a Costas sequence $a_0, a_1, \dots, a_{n-1}$, consider the circular differences 
    \begin{equation}\label{eq:circular_diff}
        a_{0+k}-a_0, \quad a_{1+k}-a_1, \quad \dots, \quad a_{n-1+k}-a_{n-1},
    \end{equation}
    where the indices $i+k$ are considered modulo $n$.
    The sequence is said to have the \textbf{shifting property} if, for each $k=1, 2, \dots, n-1$, the circular differences in \eqref{eq:circular_diff} considered modulo $n+1$ are a shift of the original sequence.
\end{definition}

\begin{prop}[\cite{moreno1992shifting}*{\S 1}]\label{prop:shifting}
    A Costas sequence has the shifting property if and only if it is Welch.
\end{prop}

Four years later, in 1996, Golomb and Moreno presented a weaker version of \Cref{conj:singly_periodic} using a condition called \emph{circular}, which lies between single periodicity and the shifting property.
That is, shifting property $\implies$ circular $\implies$ singly periodic.

\begin{definition}\label{def:circular_seq}
    A Costas sequence $a_0, a_1, \dots, a_{n-1}$ is said to be \textbf{circular} if the circular differences in \eqref{eq:circular_diff} are distinct modulo $n+1$, for each $k=1, 2, \dots, n-1$.
\end{definition}

Therefore, a circular Costas sequence is equivalent to an injective map $\varphi:\z_n \to \z_{n+1}$, with $0\not\in\im(\varphi)$, for which the difference map 
\begin{equation}\label{eq:diff_map}
    \Delta_{\varphi,k}:\z_n\to\z_{n+1}, \quad 
    i \mapsto \varphi(i+k)-\varphi(i)
\end{equation}
is injective for all $k\in\z_n\setminus\set{0}$.
Notice that this is a stronger version of the distinct difference property given in \Cref{def:costas_seq}, as here we consider the operations and equality in their corresponding group, i.e., modulo $n+1$ in the image and modulo $n$ in the domain of $\varphi$.
By making the usual identification $\z_n = \set{0,1,\dots, n-1}$, one can construct a Costas type array (with indexes starting at zero) by placing a one in position $(i,t)$ if and only if $\varphi(i) = t$.
This is an ``almost'' Costas array: it has a single one in each of the $n$ columns and a single one in each of the $n+1$ rows, except for the first row which is empty because $0\not\in\im(\varphi)$.
Costas arrays with the addition of an empty row are called \textit{extended Costas arrays} \cite{moreno2003survey}.
Extended Costas arrays equivalent to circular Costas sequences have the property that any combination of cyclic shifts in their columns and their rows will produce another extended Costas array.
In other words, repeating the pattern in both directions over the entire plane gives a doubly periodic pattern with an extended Costas array in every window of $n+1$ rows and $n$ columns.
Therefore, there is nothing special about the first row being empty.
We relax the condition $0\not\in\im(\varphi)$ and consider any injective function $\varphi:\z_n\to\z_{n+1}$ with injective difference map \eqref{eq:diff_map} to be equivalent to a circular Costas sequence, so that any row could be the (unique) empty row.

Golomb and Moreno conjectured that the only circular Costas sequences are those from the Welch construction \cite{golomb1996periodicity}*{Conjecture 1}.
\begin{conjecture}[Golomb \& Moreno \cite{golomb1996periodicity}]\label{conj:golomb_moreno}
    A Costas sequence is circular if and only if it is Welch.
\end{conjecture}

In the same 1996 article, Golomb and Moreno showed that if a Costas sequence is circular, $n+1$ is prime \cite{golomb1996periodicity}*{Theorem 1}.
This is, the sequence must have length $p-1$, where $p$ is a prime, which is consistent with the Welch sequences, thus making progress towards their conjecture.
Moreover, using the aforementioned result, any circular Costas sequence can be realized as a permutation polynomial over $\F_p$, and these are called \emph{Costas polynomials} \cite{muratovic2015characterization}.

\begin{definition}[\cite{golomb1996periodicity}*{\S III}]\label{def:costas_poly}
    A \textbf{Costas polynomial} over $\F_p$ is a permutation polynomial $f\in\F_p[x]$ with the property that $f(ax)-f(x)$ is a permutation for all $a \neq 1$.
\end{definition}

Nevertheless, a complete proof of the Golomb-Moreno conjecture remained elusive until 2015, where A. Muratovi\'c-Ribi\'c, A. Pott, D. Thompson, and Q. Wang established it by using the theory of direct product difference sets and their corresponding finite projective planes.

\begin{prop}[\cite{muratovic2015characterization}*{Theorem 3.4}]
    A Costas sequence is circular if and only if it is Welch.
\end{prop}

Most of the concepts discussed so far have been generalized to higher-dimensional arrays.
Particularly, in \cite{rubiotorres}, the present authors generalized Costas arrays to higher dimensions with a definition that includes as special cases the higher-dimensional generalizations of Costas arrays proposed by Drakakis \cite{drakakis2008higher}*{Definition 6} and by Batten and Sane \cite{batten2003permutations}*{Definition 3.2}.
A multidimensional Welch construction was proposed independently by Drakakis \cite{drakakis2008higher}*{Theorem 2} and by Moreno and Tirkel \cite{moreno2011multi}*{\S III.B} based on consecutive powers of a primitive root over an extension field.
Ortiz-Ubarri et al. \cite{ortiz2013algebraic}*{Definition 1} proposed a multidimensional analog of the circular Costas property.
Finally, Muratovi\'c-Ribi\'c et al. \cite{muratovic2015characterization}*{Definition 2.4} generalized the concept of a Costas polynomial by considering the same \Cref{def:costas_poly}, but over extension fields, not only prime fields.

In this paper, we unify all these generalizations in a coherent theoretical framework in such a way that most of the results and relations established in the two-dimensional case are shown to extend into the multidimensional setting.
In \Cref{sec:costas_circular}, we discuss a multidimensional analog of Costas arrays and relate it to circular Costas maps, a more general concept than the one proposed by Ortiz-Ubarri et al. in \cite{ortiz2013algebraic}.
In \Cref{sec:welch}, we formally define multidimensional Welch maps based on \citelist{\cite{drakakis2008higher}*{Definition 6} \cite{moreno2011multi}*{\S III.B} \cite{golomb2007status}*{Construction A}} and present three conjectures by Ortiz-Ubarri et al. on circular Costas maps.
Circular Costas maps are equivalent to abelian direct product difference sets, and this is shown in \Cref{sec:dpds}.
In \Cref{sec:costas_polys}, we present a multidimensional version of the Golomb-Moreno conjecture (\Cref{conj:golomb_moreno}), which is more general than a conjecture on Costas polynomials presented by Muratovi\'c-Ribi\'c et al. \cite{muratovic2015characterization}*{Conjecture 4.4}.
In \Cref{sec:counting}, we solve some of the conjectures by Ortiz-Ubarri et al. by counting the number of Costas polynomials in a specific class.
Finally, in \Cref{sec:shifting}, we prove a weaker version of the multidimensional Golomb-Moreno conjecture by considering a multidimensional analog of the shifting property.

\section{Multidimensional Costas arrays and circular Costas maps}
\label{sec:costas_circular}
In \cite{rubiotorres}, the authors of the present paper proposed a new definition for multidimensional Costas arrays.
For a natural number $n$, recall the notation $[n] = \set{1,\dots, n}$.

\begin{definition}\label{def:md_costas}
    Let $X = [n_1] \times \cdots \times [n_h] \subset \z^h$ and $Y = [n_{h+1}] \times \cdots \times [n_m] \subset \z^{m-h}$, for some $h$, $1 \leq h < m$.
    An \textbf{$m$-dimensional Costas array} is (equivalent to) a bijective function $\varphi:X\to Y$ with the property that, for any $k\in\z^h$, 
    \[
        \varphi(i+k)-\varphi(i) = \varphi(j+k)-\varphi(j) \implies k = 0 \textup{ or } i = j,
    \]
    for all $i,j\in X$ satisfying $i+k, j+k \in X$.
\end{definition}

Note that the restriction of \Cref{def:md_costas} to $m=2$ gives the original definition of a (two-dimensional) Costas array. As discussed in \cite{rubiotorres}, \Cref{def:md_costas} includes as special cases the multidimensional generalization of Costas arrays proposed by Drakakis \cite{drakakis2008higher}*{Definition 6} as well as the one proposed by Batten and Sane \cite{batten2003permutations}*{Definition 3.2}.
From a bijection satisfying \Cref{def:md_costas} one can construct an $m$-dimensional binary array by considering the $m$-dimensional grid defined by $X\times Y \subset \z^m$, and placing a one in the grid position $(a_1,\dots, a_h, a_{h+1}, \dots, a_m)\in X\times Y$ if and only if $\varphi(a_1,\dots,a_h) = (a_{h+1},\dots, a_m)$; the rest of the grid is filled with zeros.
This array has the property that the $m$-dimensional vectors connecting pairs of ones are all distinct.

Recall that a (two-dimensional) Costas array is equivalent to a bijection form $[n]$ to $[n]$, while a circular Costas array is equivalent to an injection from $\z_n$ to $\z_{n+1}$, both functions satisfying the distinct difference property; the former with operations over the integers, and the latter with operations over the corresponding groups.
From \Cref{def:md_costas}, an $m$-dimensional Costas array is defined by a bijection from $[n_1]\times\cdots\times[n_h]$ to $[n_{h+1}]\times\cdots\times[n_m]$, so it is natural to define an \textit{$m$-dimensional circular Costas array} by an injection 
\begin{equation}\label{eq:md_map}
    \varphi:\z_{n_1} \times \cdots \times \z_{n_h}
    \longrightarrow
    \z_{n_{h+1}} \times \cdots \times \z_{n_m}
\end{equation}
satisfying the distinct modular difference property.
From the map \eqref{eq:md_map} we can construct an $m$-dimensional binary array by filling every point $(a_1,\dots, a_h, a_{h+1}, \dots, a_m)$ on the grid $[n_1]\times\cdots\times[n_m]$ with a one if 
\[
    \varphi(a_1 \mod n_1,\ \dots,\ a_h\mod n_h) = 
    (a_{h+1}\mod n_{h+1},\ \dots, \ a_m\mod n_m),
\]
and zero otherwise.
A multidimensional array constructed from a circular Costas map preserves the Costas property periodically.
That is to say, any combination of cyclic shifts of the array entries in any direction will produce a binary array with the Costas property of distinct differences, so the vectors joining pairs of ones will be all distinct.
See \cite{thesis} for further details.

Nonetheless, any finite abelian group can be written as a direct product of cyclic groups; hence, a circular Costas map can be more generally defined as an injective function $\varphi:G_1\to G_2$, where $G_1$ and $G_2$ are finite abelian groups.
Moreover, to replicate the two-dimensional setting, we further require that $\im(\varphi) = G_2\setminus\set{g}$, for some $g\in G_2$.
This is ``almost'' a surjective function; it misses a single element in the codomain.
The missed element can be interpreted as the addition of an ``empty row'' in a multidimensional Costas arrays.
Therefore, these maps will produce multidimensional arrays which are ``almost'' Costas arrays as defined in \Cref{def:md_costas}.
Now we formally define \textit{circular Costas maps}.

From now on, $G^* = G\setminus\set{0}$.

\begin{definition}\label{def:md_circular}
    Let $G_1$ and $G_2$ be two finite abelian groups with $|G_1|+1 = |G_2|$. 
    A map $\varphi:G_1\to G_2$ is called \textbf{circular Costas} if it is injective and if, for all $k\in G_1^*$, the difference map
    \[
        \Delta_{\varphi,k}:G_1 \to G_2, \qquad i\mapsto\varphi(i+k)-\varphi(i)
    \]
    is injective.
    If $\im(\varphi) = G_2^*$, we say that $\varphi$ is a \textbf{standard} circular Costas map.
\end{definition}
\begin{definition}
    Let $\varphi_1:G_1\to G_2$ and $\varphi_2:H_1\to H_2$ be circular Costas maps.
    We say that $\varphi_1$ and $\varphi_2$ are \textbf{equivalent} if there exist two isomorphisms $\psi_1:G_1\to H_1$ and $\psi_2:G_2\to H_2$ such that
    \[
        \varphi_1(x) = y \iff \varphi_2(\psi_1(x)) = \psi_2(y).
    \]
\end{definition}

If we set $G_1 = \z_n$ and $G_2 = \z_{n+1}$, \Cref{def:md_circular} is equivalent to \Cref{def:circular_seq}.
Moreover, standard circular Costas maps with $G_1 = \z_{p^m-1}$ and $G_2 = \zm_p$ are exactly what Ortiz-Ubarri et al. defined as \textit{multidimensional periodic Costas arrays} \cite{ortiz2013algebraic}*{Definition 1}.

We want to remark that bijective maps between finite abelian groups satisfying the distinct difference property do not exist \citelist{\cite{drakakis2009apn}*{Theorem 1} \cite{rubiotorres}*{Theorem 2}}.
Hence, circular Costas maps are the closest we can get to a bijection between finite abelian groups with the distinct difference property.

From a single circular Costas map $\varphi:G_1\to G_2$, it is possible to construct several multidimensional arrays, each in different dimensions, depending on the number of factors in the chosen decomposition of $G_1$ and $G_2$ as product of cyclic groups.
This is shown in the next example.

\begin{example}\label{ex:decomp}
    Let $\alpha$ be a primitive element in $\F_{25}$.
    Consider the map 
    \[
        \varphi:\z_{24}\to(\F_{25},+), \qquad i\mapsto\alpha^i.
    \]
    This is a circular Costas map (shown in \Cref{prop:md_welch} below).
    Let $\psi_1:\z_{24}\to\z_8\times\z_3$ and $\psi_2:\F_{25}\to\z_5\times\z_5$ be group isomorphisms.
    By making the identification $\z_n=\set{0,1,\dots, n-1}\subset\z$, the map 
    \[
        \varphi_1:\z_{24}\to\z_5\times\z_5, \qquad i\mapsto\psi_2(\alpha^i),
    \]
    defines a three-dimensional binary array of size $24 \times 5 \times 5$ preserving the Costas property periodically in all three directions.
    
    On the other hand, the map
    \[
        \varphi_2:\z_8\times\z_3\to\z_5\times\z_5, \qquad \psi_1(i)\mapsto\psi_2(\alpha^i), \ i\in\z_{24},
    \]
    defines a 4-dimensional binary array of size $8\times 3 \times 5 \times 5$ preserving the Costas property periodically in all four directions.
\end{example}

\section{Multidimensional Welch maps}\label{sec:welch}
Consider a primitive element $\alpha$ of $\F_q$, the finite field with $q = p^m$ elements, $p$ prime.
Write the elements of $\F_{q}$ as $m$-tuples over $\z_p$ with respect to $\alpha$, so that 
\[
    \alpha^i = c_0 + c_1\alpha + \cdots + c_{m-1}\alpha^{m-1} \mapsto (c_0, c_1, \dots, c_{m-1})\in\z_p^m. 
\]
By making the usual identification $\z_n = \set{0,1,\dots, n-1} \subset \z$, for fixed $c\in\z$, in the $(m+1)$-dimensional grid of size $q-1\times p \times\cdots\times p$, place a one in position $(i,\alpha^{i+c})$ for all $i\in\set{0,1,\dots,q-2}$, and zeros elsewhere.
The resulting array satisfies the Costas property of distinct differences.
This construction was described independently by Drakakis \cite{drakakis2008higher}*{Theorem 2} and by Moreno and Tirkel \cite{moreno2011multi}*{Secton III.B}.
Moreover, let $L(x)$ be a linearized permutation polynomial of $\F_q$, so it is a permutation polynomial of the form $L(x) = \sum_{j=0}^{m-1}c_jx^{p^j}$.
Construct a similar array as above but where the ones are placed in positions $(i, L(\alpha^i))$.
Again, we consider $L(\alpha^i)$ as an $m$-tuple in base $\alpha$.
The latter $(m+1)$-dimensional array defined by the mapping $i \mapsto L(\alpha^{i+c})$ preserves the Costas property of distinct differences, as shown by Golomb and Gong in \cite{golomb2007status}*{Theorem 1}.
Based on the remark in \cite{golomb2007status}*{p.~4264}, we call the map $\varphi(i) = L(\alpha^{i+c})$ a \textit{multidimensional Welch map}.

\begin{definition}\label{def:md_welch}
    Let $\alpha$ be a primitive element in $\F_q$, with $q$ a prime power, let $L\in\F_q[x]$ be linearized permutation polynomial, and let $c\in\z$.
    A \textbf{multidimensional Welch map} is any map equivalent to \begin{equation}\label{eqwelch}
        \varphi:\z_{q-1}\to(\F_q,+), \qquad i\mapsto L(\alpha^{i+c}).
    \end{equation}
\end{definition}

\begin{prop}\label{prop:md_welch}
    Multidimensional Welch maps are standard circular Costas.
\end{prop}
\begin{proof}
    Let $\varphi:G_1\to G_2$ be a multidimensional Welch map.
    Without loss of generality, we may assume $\varphi$ is as in \eqref{eqwelch}.
    A linearized permutation polynomial is, by definition, a bijective linear operator from $\F_q$ to itself.
    Since $\alpha\in\F_q$ is primitive and $L$ is a linearized permutation polynomial, $\varphi$ is injective and $0\not\in\im(\varphi)$.
    For fixed $k\in\z_{q-1}^*$, consider the difference map
    \[
        \Delta_{\varphi,k}:\z_{q-1}\to\F_q, \qquad i\mapsto\varphi(i+k)-\varphi(i).
    \]
    Then, by considering $i,j$, and $k$ as their canonical integer representatives (notice $k\neq 0$),
    \begin{align*}
    \Delta_{\varphi,k}(i) = \Delta_{\varphi,k}(j)
        &\implies \varphi(i+k) - \varphi(i) = \varphi(j+k) - \varphi(j) \\
        &\implies L(\alpha^{i+c+k}) - L(\alpha^{i+c}) = L(\alpha^{j+c+k})-L(\alpha^{j+c})\\
        &\implies L(\alpha^{i+c+k}-\alpha^{i+c}) = L(\alpha^{j+c+k}-\alpha^{j+c})\\
        &\implies \alpha^i\alpha^c(\alpha^k-1) = \alpha^j\alpha^c(\alpha^k-1)\\
        &\implies \alpha^i = \alpha^j\\
        &\implies i=j.
    \end{align*}
    Hence, the difference map is injective.
    We conclude that $\varphi$ is a standard circular Costas map.
\end{proof}

\begin{remark}
    When $q$ is a prime, $q=p$, the only monic linearized permutation polynomial $L\in\F_p$ is $L(x) = x$, which is the identity function.
    In this case, the map $\varphi$ in \eqref{eqwelch} has the rule $i\mapsto L(\alpha^{i+c}) = \alpha^{i+c}$, and we get back the two-dimensional Welch construction (\Cref{prop:welch}).
    If $L$ is not monic, the non-zero leading coefficient can be expressed as a power of $\alpha$, so it will be absorbed by the constant $c$.
\end{remark}

In \cite{ortiz2013algebraic}, Ortiz-Ubarri et al. gave four conjectures regarding circular Costas maps, three of which we state here.
The other conjecture is out of the scope of this paper.

\begin{conjecture}[\cite{ortiz2013algebraic}*{Conjecture 4}]\label{conj:cheo4}
    If $\varphi:G_1\to G_2$ is circular Costas, $G_2$ is elementary abelian.
\end{conjecture}
\begin{conjecture}[\cite{ortiz2013algebraic}*{Conjecture 2}]\label{conj:cheo2}
    The number of standard circular Costas maps $\varphi:\z_{p^2-1}\to\z_p^2$ is at least
    \begin{equation}\label{eq:conj2}
        \frac{\phi(p^2-1)}{2}(p^2-p)(p^2-1),
    \end{equation}
    where $\phi$ is the Euler's totient function.
\end{conjecture}
\begin{conjecture}[\cite{ortiz2013algebraic}*{Conjecture 3}]\label{conj:cheo3}
    For $m\geqslant 3$, the number of standard circular Costas maps $\varphi:\z_{p^m-1}\to\zm_p$ is at least
    \begin{equation}\label{eq:conj3}
        \frac{\phi(p^m-1)}{m}(\phi(p^m-1)-1)(p^m-1)(p-1)^{m-1}p^{m-1},
    \end{equation}
    where $\phi$ is the Euler's totient function.
\end{conjecture}

Conjectures \ref{conj:cheo4} and \ref{conj:cheo2} will be proved in Sections \ref{sec:dpds} and \ref{sec:counting}, respectively.
We will show that \Cref{conj:cheo3} is true for $m\geqslant 4$, but we think it is false for $m=3$, and we discuss why in \Cref{sec:counting}.

\section{Circular Costas maps and direct product difference sets}\label{sec:dpds}
In this section we show that circular Costas maps are equivalent to abelian direct product difference sets.
This, together with a result by Jungnickel and de Resmini, is used to prove Conjecture 4 in \cite{ortiz2013algebraic} (here, \Cref{conj:cheo4}).
We want to remark that the connection between circular Costas arrays and direct product difference sets is not new, as seen in \citelist{\cite{drakakis2010nonlinearity} \cite{muratovic2015characterization}}.

\begin{definition}[\cite{ganley1977direct}*{p.~321}]\label{def:dp_difference_set}
    Let $G = A\times B$ be the direct product of two groups $A$ and $B$ (written additively), with $|A|=n-1$ and $|B|=n$, where $n\geqslant 3$.
    Let $D\subseteq G$ be such that every element of $G\setminus\set{(A\times\set{0}) \cup (\set{0}\times B)}$ can be uniquely represented in the form $d_j-d_i$, where $d_i, d_j \in D$.
    Furthermore, suppose that no nonidentity element of $(A\times\set{0}) \cup (\set{0}\times B)$ can be so represented.
    We call $D$ a \textbf{direct product difference set} in $G$ of order $n$.
\end{definition}

Notice that, in a direct product difference set $D$ in the group $G = A\times B$, if $(a_1,b_1)\in D$ and $(a_2,b_2)\in D$ are distinct, then $a_1\neq a_2$ and $b_1\neq b_2$.
Otherwise, an element of $(A\times\set{0}) \cup (\set{0}\times B)$ could be represented as the difference of elements in $D$.
By the previous observation, a simple count gives that, if $D\subseteq G$ is a direct product difference set of order $n$, $|D| = n-1$.

\begin{example}
    Let $G = \z_4\times\z_5$. 
    Consider the subset
    \[
        D = \set{(0,2), \ (1,4), \ (2,3), \ (3,1)} \subseteq G.
    \]
    Take all the differences between distinct elements of $D$:
    \begin{center}
    \begin{tabular}{cc}
        $(0,2)-(1,4)=(3,3)$\qquad\ & \ \qquad$(1,4)-(0,2)=(1,2)$ \\
        $(0,2)-(2,3)=(2,4)$\qquad\ & \ \qquad$(1,4)-(2,3)=(3,1)$ \\
        $(0,2)-(3,1)=(1,1)$\qquad\ & \ \qquad$(1,4)-(3,1)=(2,3)$ \\
        $(2,3)-(0,2)=(2,1)$\qquad\ & \ \qquad$(3,1)-(0,2)=(3,4)$ \\
        $(2,3)-(1,4)=(1,4)$\qquad\ & \ \qquad$(3,1)-(1,4)=(2,2)$ \\
        $(2,3)-(3,1)=(3,2)$\qquad\ & \ \qquad$(3,1)-(2,3)=(1,3)$. 
    \end{tabular}
    \end{center}
    Every element of $\z_4^*\times\z_5^*$ appears exactly once and no element of $(\z_4\times\set{0})\cup(\set{0}\times\z_5)$ appears as a difference from elements in $D$; therefore, $D$ is a direct product difference set in $G$ of order $5$.
\end{example}

Jungnickel and de Resmini \cite{jungnickel2002another} proved this remarkable result about direct product difference sets.

\begin{prop}[\cite{jungnickel2002another}*{Theorem 1}]\label{prop:prime_power}
	Let $D$ be an order $n$ direct product difference set in the abelian group $G = A \times B$, where $|A| = n-1$ and $|B|=n$.
	The order $n$ must be a prime power and the group $B$ is elementary abelian.\footnote{
		Theorem 1 in \cite{jungnickel2002another} is not stated in terms of direct product difference sets, but rather in terms of collineation groups of finite projective planes.
	}
\end{prop}

Two direct product difference sets $D_1$ and $D_2$ in the abelian group $G = G_1\times G_2$ are said to be \textit{equivalent} if there is an element $(a,b)\in G$ and a group automorphism $\psi\in\aut(G)$, such that $D_2 = \set{(a+x,b+y) \st (x,y)\in\psi(D_1)} $.
From a direct product difference set $D\subseteq G$, define a function $\varphi:G_1 \to G_2$ by setting
\begin{equation}\label{eq:associated_func}
    \varphi(i) = j \iff (i,j)\in D.
\end{equation}
A function defined as in \eqref{eq:associated_func} is said to be \textit{associated with the direct product difference set} $D$.

\begin{theorem}\label{thm:md_circular_dp_difference}
    $D$ is a direct product difference set of order $n = |G_2| = |G_1|+1$ in the abelian group $G_1\times G_2$ if and only if the associated function $\varphi:G_1\to G_2$ is a circular Costas map.
\end{theorem}
\begin{proof}
    Let $D = \set{(i,\varphi(i)) \st i\in G_1}$ be a direct product difference set in $G_1 \times G_2$, for some map $\varphi: G_1 \to G_2$.
    The map $\varphi$ must be injective; otherwise, one could write a nonidentity element on the subgroup $G_1\times\set{0}$ as a difference from elements in $D$.
    Since $D$ is a direct product difference set, for any $k\neq 0$, the difference $\varphi(i+k)-\varphi(i)$ must run through $G_2^*$ as $i$ runs through $G_1$.
    Hence, the difference map $\Delta_{\varphi,k}$ is injective for all $k\neq 0$.
    We conclude that $\varphi$ is circular Costas.
    
    Conversely, let $\varphi:G_1\to G_2$ be a circular Costas map, and let 
    \[
        D = \set{(i,\varphi(i)) \st i\in G_1}.
    \]
    Take $(k,h)\in G_1 \times G_2$.
    Assume $(k,h)\not\in (G_1\times\set{0}) \cup (\set{0} \times G_2)$.
    Since $h\neq 0$ and $k\neq 0$, the injectivity of $\Delta_{\varphi,k}$ and the fact that $\im(\Delta_{\varphi,k}) = G_2^*$ ensures a unique $i\in G_1$ such that $h = \varphi(i+k)-\varphi(i)$.
    Then $(k,h) = (i+k,\varphi(i+k))-(i,\varphi(i))$.
    Therefore, $(k,h)$ can be written uniquely as a difference from elements in $D$.
    If $(j-i,\varphi(j)-\varphi(i)) = (0,h)$ then $i = j$ and $h = 0$.
    On the other hand, if $(j-i,\varphi(j)-\varphi(i)) = (k,0)$, then $i = j$ because $\varphi$ is injective.
    We conclude that $D$ is a direct product difference set in $G_1\times G_2$.
\end{proof}

The above theorem is similar to Lemma 3.3 in \cite{muratovic2015characterization}.
However, our statement is more general because theirs is in terms of Costas polynomials.
In \Cref{sec:costas_polys}, we will see that Costas polynomials are a special types of circular Costas maps.

The next theorem establishes affirmatively \Cref{conj:cheo4} (Conjecture 4 in \cite{ortiz2013algebraic}).

\begin{theorem}[\Cref{conj:cheo4}]\label{thm:conj4}
    If $\varphi:G_1 \to G_2$ is circular Costas, $G_2$ is elementary abelian.
    That is, $G_2 \cong \zm_p$ for some natural $m$ and prime $p$.
\end{theorem}
\begin{proof}
    Let $\varphi:G_1 \to G_2$ be circular Costas, and let $|G_1|+1 = |G_2| = n$. By \Cref{thm:md_circular_dp_difference}, $D = \set{(i,\varphi(i)) \st i\in G_1}$ is a direct product difference set in $G_1 \times G_2$ of order $n$.
    \Cref{prop:prime_power} shows that $G_2$ must be elementary abelian.
\end{proof}

In light of \Cref{thm:conj4}, we give the following definition.

\begin{definition}
    A circular Costas map $\varphi:G_1\to G_2$ has \textbf{order} $q$ and \textbf{degree} $m$ if $G_2 \cong\zm_p$.
\end{definition}

\begin{cor}
    Any circular Costas map has prime power order and unique degree.
\end{cor}

\section{Costas polynomials}\label{sec:costas_polys}

In the context of \Cref{thm:conj4}, any circular Costas map with cyclic domain is equivalent to an injective map $f:(\F_q^*,\times)\to(\F_q,+)$, where now the domain uses multiplicative notation.
If the circular Costas map $f$ is standard, $f(x) \neq 0$ for all $x\in\F_q^*$.
We can extend the domain and the image of the map $f$ by setting $f(0) = 0$.
Such $f$ is a bijection from $\F_q$ to itself, which can be expressed, by Lagrange interpolation, as a unique polynomial function of degree at most $q-1$.
Hence, the concept of Costas polynomials over extension fields also makes sense in the context of multidimensional circular Costas maps.

\begin{definition}[\cite{muratovic2015characterization}*{Definition 2.4}]\label{def:md_costas_poly}
    Let $q$ be a prime power.
    A polynomial $f\in\F_q[x]$ is a \textbf{Costas polynomial} if $f(0) = 0$ and $f(dx)-f(x)$ is a permutation for all $d\in\F_q$, $d\neq 1$.
\end{definition}

Muratovi\'c-Ribi\'c et al. mentioned that ``while we have defined Costas polynomials over any finite field, we emphasize that they are equivalent to circular Costas sequences only over prime fields.'' \cite{muratovic2015characterization}*{p.~320}
However, with our definition of circular Costas maps and the next lemma, we extend the equivalence to include Costas polynomials over extension fields.

\begin{lemma}\label{lemma:one-to-one_correspondence}
    Let $\varphi:G_1\to G_2$ be a standard circular Costas map of order $q$.
    Then $\varphi$ is equivalent to a Costas polynomial $f\in\F_q[x]$ if and only if $G_1$ is cyclic.
    Moreover, there is a one-to-one correspondence between the set of order $q$ standard circular Costas maps with cyclic domain and the set of Costas polynomials of degree at most $q-1$. 
\end{lemma}
\begin{proof}
    Let $\varphi:G_1\to G_2$ be a standard circular Costas map of order $q$.
    First assume $\varphi:G_1\to G_2$ is equivalent to a Costas polynomial $f\in\F_q[x]$.
    By the definition of a Costas polynomial, if we restrict the domain of $f$ to the multiplicative group of $\F_q$, then $f:(\F_q^*,\times)\to(\F_q,+)$ is a standard circular Costas map (with the domain written multiplicatively).
    Therefore, $G_1\cong (\F_q^*,\times)$ and $G_1$ is cyclic.
    
    Conversely, assume $G_1$ is cyclic.
    By \Cref{thm:conj4}, $G_2 \cong (\F_q,+)$ for some prime power $q$, and by the definition of a circular Costas map, $|G_1| = q-1$.
    Since $G_1$ is cyclic, we have $G_1 \cong (\F_q^*,\times)$.
    Let $\psi_1:G_1\to(\F_q^*,\times)$ and $\psi_2:G_2\to(\F_q,+)$ be two group isomorphisms.
    Define the map $f:\F_q^*\to\F_q$, by the rule
    \begin{equation}\label{eq:rule}
        \varphi(i) = j \iff f(\psi_1(i)) = \psi_2(j).
    \end{equation}
    Notice that $f$ is injective and $0\not\in\im(f)$ because $\varphi$ is standard.
    Extend the domain of $f$ from $\F_q^*$ to $\F_q$ by setting $f(0)=0$.
    Since $\varphi$ is circular Costas, $f$ with the rule \eqref{eq:rule} is a Costas polynomial.
    We conclude that $\varphi$ is equivalent to the Costas polynomial $f$.
    
    Notice that, by Lagrange interpolation, for fixed isomorphisms $\psi_1$ and $\psi_2$, $f$ is unique if we restrict it to have degree at most $q-1$.
    The isomorphisms ensure a one-to-one correspondence between order $q$ standard circular Costas maps with cyclic domain and Costas polynomials of degree at most $q-1$.
\end{proof}

A direct product difference set of order $n$ in a group $G_1\times G_2$ is equivalent to a projective plane of order $n$ admitting a quasiregular collineation group isomorphic to $G_1 \times G_2$ \cite{ganley1977direct}*{Theorem 4.3}.
It is conjectured that projective planes that are equivalent to abelian direct product difference sets are desarguesian \cite{jungnickel2002another}*{p.~216}.
By \cite{pott1994projective}*{Theorem 2.3}, the following is an equivalent conjecture in terms of circular Costas maps.

\begin{conjecture}\label{conj:costas_polys}
    Every standard circular Costas map is equivalent to a Costas polynomial $f\in\F_q$ with $f(x) \equiv L(x^s) \pmod{x^q-x}$, where $L(x)\in\F_q[x]$ is a linearized permutation polynomial and $\gcd(s,q-1) = 1$.
\end{conjecture}

\begin{remark}
    A polynomial $f$ as in \Cref{conj:costas_polys} is indeed a Costas polynomial; shown in \cite{muratovic2015characterization}*{Theorem 4.1}.
\end{remark}

\Cref{conj:costas_polys} is twofold.
The first part is the assertion that every standard circular Costas map $\varphi$ is equivalent to a Costas polynomial; by \Cref{lemma:one-to-one_correspondence}, that is the same as saying that the domain of $\varphi$ is cyclic.
The second part of \Cref{conj:costas_polys} regards the form of the Costas polynomials and is the same as \cite{muratovic2015characterization}*{Conjecture 4.4}.
We claim that the polynomials in \Cref{conj:costas_polys} correspond to the multidimensional Welch maps of \Cref{def:md_welch}.

Take a primitive element $\alpha\in\F_q$, a linearized permutation polynomial $L(x)\in\F_q[x]$, a constant $c\in\z$, and consider the multidimensional Welch map
\[
    \varphi:\z_{q-1}\to(\F_q,+), \qquad i\mapsto L(\alpha^{i+c}).
\]
To obtain the Costas polynomial corresponding to $\varphi$, we need a group isomorphism $\psi:\z_{q-1}\to\F_q^*$.
Take $\psi(i) = \beta^i$, where $\beta$ is a primitive element in $\F_q$.
The Costas polynomial $f\in\F_q[x]$ corresponding to the multidimensional Welch map is the polynomial interpolating the points $(\psi(i), \varphi(i))$ for $i\in\z_{q-1}$.
Therefore,
\[
    f(\psi(i)) = \varphi(i) \implies f(\beta^i) = L(\alpha^{i+c}).
\]
But $\alpha = \beta^s$, for some $s$ with $\gcd(s,q-1)=1$.
Hence, by substituting $\alpha=\beta^s$ and taking $a = \alpha^c$, $f(x) = L(ax^s)$.
Moreover, we can choose a suitable linearized permutation polynomial $L'$ such that $L(ax) = L'(x)$. 
We end up with $f(x) = L'(x^s)$, where $L'$ is a linearized permutation polynomial and $\gcd(s,q-1) = 1$.
Therefore, the polynomials in \Cref{conj:costas_polys} correspond to multidimensional Welch maps.
\Cref{conj:costas_polys} may be rephrased to assert that any circular Costas map is multidimensional Welch; a multidimensional version of the Golomb-Moreno conjecture (\Cref{conj:golomb_moreno}).

\begin{conjecture}[Equivalent to \Cref{conj:costas_polys}]\label{conj:md_circular_implies_Welch}
    All standard circular Costas maps are multidimensional Welch.
\end{conjecture}

\section{Counting Costas polynomials}\label{sec:counting}

By \Cref{lemma:one-to-one_correspondence}, counting the number of Costas polynomials in $\F_q[x]$ of degree at most $q-1$ is the same as counting the number of standard circular Costas maps with cyclic domain and order $q$.
Of course, if \Cref{conj:costas_polys} were to be true, the number of Costas polynomials (of degree at most $q-1$) and hence the number of standard circular Costas maps with cyclic domain should be exactly the number of inequivalent Costas polynomials of the type $L(x^s)$, for $L(x)\in\F_q[x]$ a linearized permutation polynomial and $\gcd(s,q-1) = 1$.
However, regardless the veracity of \Cref{conj:costas_polys}, by counting the number of Costas polynomials of the type $L(x^s)$, we will obtain a lower bound on the number of standard circular Costas maps with cyclic domain and, therefore, a lower bound on general standard circular Costas maps.

In the following, $\phi$ denotes Euler's totient function.

\begin{theorem}\label{thm:number_of_costas}
    The number of distinct Costas polynomials over $\F_{p^m}$ (distinct modulo $x^{p^m}-x$) of the form $L(x^s)$, for some linearized permutation polynomial $L(x)\in\F_{p^m}[x]$ and $\gcd(s,p^m-1) = 1$, is
    \begin{equation}\label{eq:number_of_costas}
        \frac{\phi(p^m-1)}{m}\prod_{k=0}^{m-1} (p^m-p^k).
    \end{equation}
\end{theorem}
\begin{proof}
Let $q = p^m$, $p$ prime.
The polynomials of the form $\sum_{i=0}^{m-1}c_ix^{p^i} \in \F_q[x]$ are linearized polynomials.
From \cite{zhou2008remark}*{Theorem 2.1}, there are 
\[
    \prod_{k=0}^{m-1} (p^m-p^k)
\]
distinct (modulo $x^q-x$) linearized permutation polynomials over $\F_q$.

Assume $L(x) = \sum_{i=0}^{m-1}c_ix^{p^i}$ is a linearized permutation polynomial.
Since we will consider the polynomials modulo $x^q-x$, consider only $s < q-1$ with $\gcd(s,q-1) = 1$.
There are $\phi(q-1)$ choices of $s$ for which $L(x^s)$ is a Costas polynomial.

Now we have to count the number of distinct ways we can produce the same Costas polynomial modulo $x^q-x$.
Let $L_0(x) = \sum_{i=0}^{m-1}c_ix^{p^i}$ be a linearized permutation polynomial over $\F_q$ and fix $s_0$ with $\gcd(s_0,q-1)=1$.
Assume there is a linearized permutation polynomial $L(x) = \sum_{i=0}^{m-1}b_ix^{p^i}$ and a number $s<q-1$ with $\gcd(s, q-1)=1$, such that
\[
    L(x^s) \equiv L_0(x^{s_0}) \pmod{x^q-x}.
\]
By comparing the powers of $x$, we must have $s \equiv s_0p^k \pmod{q-1}$ for some integer $k$ satisfying $0 \leqslant k \leqslant m-1$.
Then,
\[
    L(x^s) \equiv L(x^{s_0p^k}) \equiv \sum_{i=0}^{m-1} b_ix^{s_0p^{k+i}} \equiv L_0(x^{s_0}) \pmod{x^q-x}.
\]
Therefore, by comparing coefficients, $b_i = c_{i+k}$, with the index $i+k$ considered modulo $m$, implying $L(x^{p^k}) \equiv L_0(x)$.
The converse is clearly true.

This shows that $L(x^s) \equiv L_0(x^{s_0})$ if and only if $L(x^{p^k}) \equiv L_0(x)$ and $s\equiv s_0p^k \pmod{q-1}$, for some $k$, where $0\leqslant k \leqslant m-1$.
We conclude that there are exactly $m$ distinct ways to obtain the same Costas polynomial, one for each selection of $k$. 
\Cref{eq:number_of_costas} follows from the multiplicative principle.
\end{proof}

\Cref{thm:number_of_costas} above makes all the heavy lifting required to prove \Cref{conj:cheo2}.

\begin{cor}[\Cref{conj:cheo2}]\label{cor:cheo2}
    The number of standard circular Costas maps $\varphi:\z_{p^2-1}\to\z_p^2$ is at least
    \[
        \frac{\phi(p^2-1)}{2}(p^2-p)(p^2-1).
    \]
\end{cor}
\begin{proof}
Follows from \Cref{lemma:one-to-one_correspondence} and \Cref{thm:number_of_costas} with $m=2$.
\end{proof}

For $m\geqslant 3$, let us compare the lower bound on the number of standard circular Costas maps given in \Cref{thm:number_of_costas} with the lower bound given in \Cref{conj:cheo3}.
Consider the ratio between the two lower bounds:
\begin{equation}\label{eq:ratio}
    R(p,m) = \frac{\frac{\phi(p^m-1)}{m}\prod_{k=0}^{m-1} (p^m-p^k)}{\frac{\phi(p^m-1)}{m}(\phi(p^m-1)-1)(p^m-1)(p-1)^{m-1}p^{m-1}}.
\end{equation}
If $R(p,m) \geqslant 1$ for all prime $p$ and all integer $m\geqslant 3$, \Cref{conj:cheo3} would be proved.
Let us simplify the expression in \eqref{eq:ratio}.
For $m\geqslant 3$, we can rewrite
\begin{align}
    \prod_{k=0}^{m-1} (p^m-p^k) &= (p^m-1)\prod_{k=1}^{m-1}p^k(p^{m-k}-1)\nonumber\\
        &=(p^m-1)p^{\sum_{k=1}^{m-1}k}\prod_{k=1}^{m-1}(p^{m-k}-1)\nonumber\\
        &=(p^m-1)p^{\frac{m(m-1)}{2}}(p-1)^{m-1}\prod_{k=1}^{m-1}\frac{p^{m-k}-1}{p-1}.\label{eq3}
\end{align}
Moreover, the product on the far right of \eqref{eq3} can be bounded below by
\begin{equation}\label{eq4}
    \prod_{k=1}^{m-1}\frac{p^{m-k}-1}{p-1} = (p+1)\prod_{\substack{k=1 \\ k\neq m-2}}^{m-1}\frac{p^{m-k}-1}{p-1} \geqslant p+1,
\end{equation}
where equality holds only when $m=3$.
Therefore, by using \eqref{eq3}, we can simplify \eqref{eq:ratio} to obtain
\begin{equation}\label{eq5}
    R(p,m) = \frac{p^{\frac{(m-1)(m-2)}{2}}\prod_{k=1}^{m-1}\frac{p^{m-k}-1}{p-1}}{\phi(p^m-1)-1}.
\end{equation}

From \eqref{eq4} and \eqref{eq5}, for $m\geqslant 4$,
\begin{align*}
    R(p,m) > \frac{p^{\frac{(m-1)(m-2)}{2}}(p+1)}{\phi(p^m-1)-1} > \frac{p^m}{\phi(p^m-1)-1} > 1.
\end{align*}
This shows that, for $m\geqslant 4$, the lower bound given in \Cref{thm:number_of_costas} for the number of standard circular Costas maps of degree $m$ is better than the one given in \Cref{conj:cheo3}.
Hence, we have shown \Cref{conj:cheo3}, except for $m=3$.

\begin{prop}\label{prop:cheo3}
    For $m\geqslant 4$, the number of standard circular Costas maps $\varphi:\z_{p^m-1}\to\zm_p$ is greater than
    \begin{equation*}
        \frac{\phi(p^m-1)}{m}(\phi(p^m-1)-1)(p^m-1)(p-1)^{m-1}p^{m-1},
    \end{equation*}
    where $\phi$ is the Euler's totient function.
\end{prop}

\begin{figure}[ht]
    \centering
    \includegraphics[scale=.2]{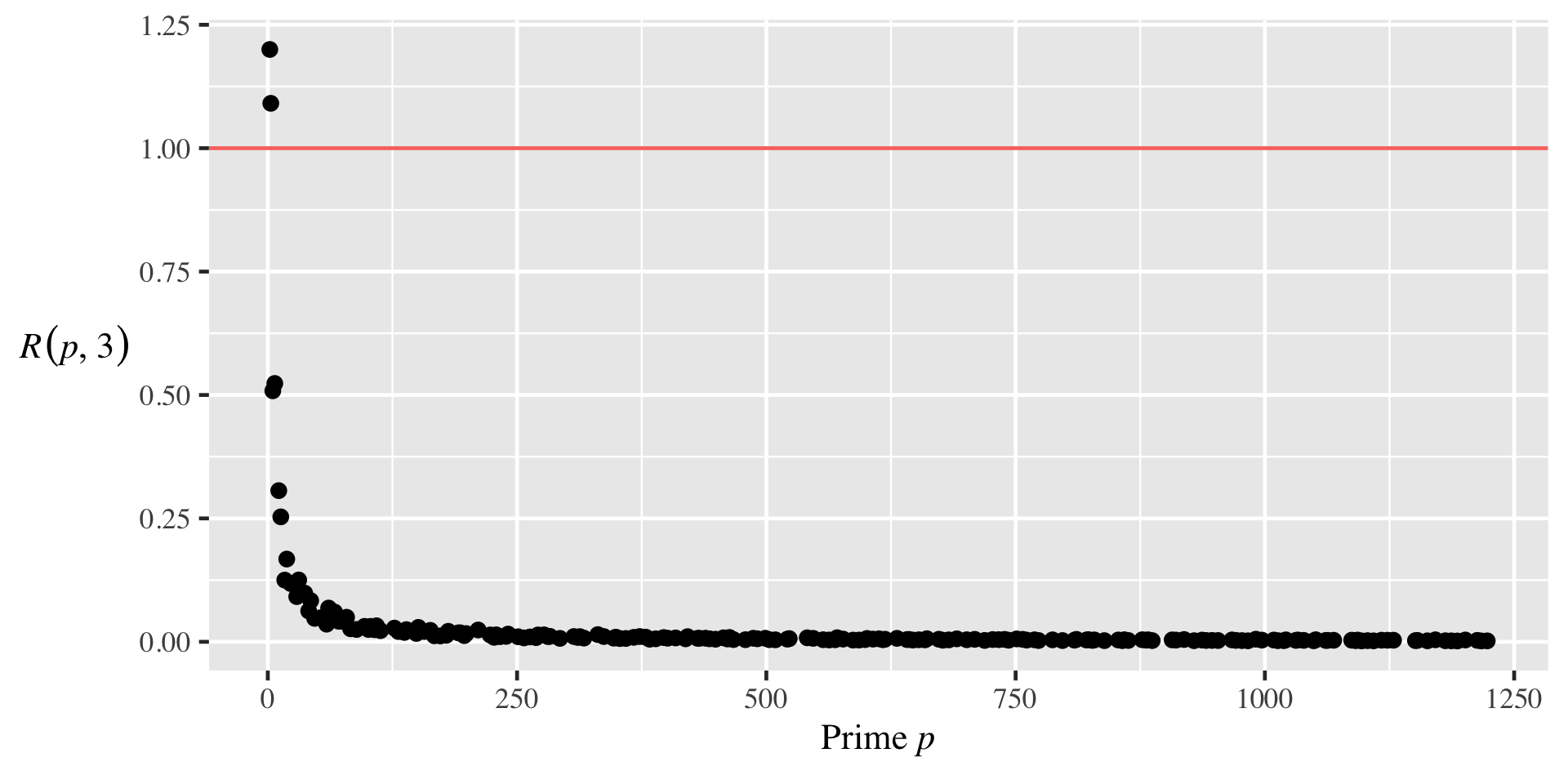}
    \caption{Asymptotic behavior of $R(p,3)$}
    \label{fig:plot}
\end{figure}

The case $m=3$ yields a different conclusion.
Consider 
\[
    R(p,3) = \frac{p(p+1)}{\phi(p^3-1)-1}.
\]
It happens that $R(2,3) = 6/5 > 1$ and $R(3,3) = 12/11 > 1$.
However, $R(5,3) = 30/59 < 1$.
In fact, it appears that $R(p,3) < 1$ for $p\geqslant 5$ (see \Cref{fig:plot}).
Thus, for $m=3$, the lower bound given in \Cref{thm:number_of_costas} is almost never better than the one given in \Cref{conj:cheo3}.

If \Cref{conj:cheo3} ought to be true for $m=3$, it would imply that \Cref{conj:costas_polys} is false.
The perception of the authors is that \Cref{conj:cheo3} is false (for $m=3$).
The falsehood of \Cref{conj:cheo3} could be checked computationally by constructing all the standard circular Costas maps for the smallest $p$ for which $R(p,3) < 1$, i.e., $p=5$.
However, the computer would need to check the Costas property for all the $(5^3-1)! = 124!$ possible permutation polynomials of $\F_{5^3}$ fixing zero.
We do not have the computational resources for this computation to finish in a reasonable time.

\section{Shifting Costas polynomials}\label{sec:shifting}
In this section we will prove a weaker version of \Cref{conj:costas_polys}.
We do it by following O. Moreno's idea in the context of two-dimensional circular Costas arrays: we consider a multidimensional extension of the shifting property (\Cref{def:moreno_shifting}).

\begin{definition}
    Let $q$ be a prime power.
    A \textbf{shifting Costas polynomial} is a permutation polynomial $f\in\F_q[x]$ with the property that, for any $d\neq1$, there exist $a\in\F_q^*$ such that $f(dx)-f(x) = f(ax)$.
\end{definition}

From a Costas polynomial $f\in\F_q[x]$, choose a primitive element $\alpha\in\F_q$ and consider the sequence
\[
    s_0 = f(\alpha^0), \quad s_1 = f(\alpha^1), \quad \dots, \quad s_{q-2} = f(\alpha^{q-2}).
\]
Let $d = \alpha^c$ and $x=\alpha^i$.
Then $f(dx)-f(x) = s_{i+c}-s_i$, where the indices are taken modulo $q-1$.
Notice that, if $f$ is a shifting Costas polynomial and $a = \alpha^k$,
\[
    s_{i+c}-s_i = s_{i+k}.
\]
Hence, the sequence of differences is a circular shift of the original sequence. 
This is why we call these polynomials ``shifting'', as they produce a multidimensional generalization of the shifting property given in \Cref{def:moreno_shifting}, but expressed in terms of their Costas polynomial.

\begin{lemma}\label{lemma:shifting_is_costas}
    A shifting Costas polynomial is a Costas polynomial.
\end{lemma}
\begin{proof}
    Let $f\in\F_q[x]$ be a shifting Costas polynomial.
    If $a\in\F_q^*$, $x\mapsto ax$ is a permutation of $\F_q$.
    The polynomial $f(ax)$ is the composition of permutations; hence, it is also a permutation of $\F_q$.
    Then, $f(dx)-f(x) = f(ax)$ is a permutation for all $d\in\F_q, d\neq 1$.
    Moreover, $f(ax)=f(dx)-f(x)$ has no constant term, implying $f(0) = 0$.
    We conclude that $f$ is a Costas polynomial.
\end{proof}

\begin{lemma}\label{lemma:composition}
    Let $f \in \F_q[x]$ be a shifting Costas polynomial.
    The polynomial $f(x^s)$ is a shifting Costas polynomial if and only if $\gcd(s,q-1)=1$. 
\end{lemma}
\begin{proof}
    If $\gcd(s,q-1)\neq 1$, then $x^s$ is not a permutation monomial. 
    In particular, $x^s$ is not a surjective function over $\F_q$, so $f(x^s)$ is not a permutation.
    Hence, $f(x^s)$ is not a shifting Costas polynomial.
    
    Now assume $s$ is a positive integer with $\gcd(s,q-1) = 1$ and fix $d\in\F_q$, $d\neq 1$.
    Define $g(x) = f(x^s)$.
    We have $g(0) = f(0) = 0$.
    Moreover, $x\mapsto x^s$ is a permutation monomial, and since $f$ is a permutation by assumption, $g(x)$ is the composition of permutations, hence a permutation.
    
    Notice that $d^s\neq 1$.
    Since $f$ is shifting Costas, there exists $a\in\F_q^*$ for which
    \[
        f(d^sx)-f(x) = f(ax).
    \]
    Since $\gcd(s,q-1)=1$, there is some $a_0\in\F_q$ such that $a_0^s = a$.
    Therefore,
    \[
        g(dx)-g(x) = f(d^sx^s)-f(x^s) = f(ax^s) = g(a_0x).
    \]
    We conclude that $g(x) = f(x^s)$ is a shifting Costas polynomial.
\end{proof}

And now a weaker version of \Cref{conj:costas_polys}, which is also a multidimensional extension of \Cref{prop:shifting}.

\begin{theorem}\label{thm:md_shifting}
    Let $q$ be a prime power.
    A polynomial $f \in \F_q[x]$ is a shifting Costas polynomial if and only if $f(x) \equiv L(x^s) \pmod{x^q-x}$, where $L(x)\in\F_q[x]$ is a linearized permutation polynomial and $\gcd(s,q-1) = 1$.
\end{theorem}
\begin{proof}
    Let $f\in\F_q$ be a shifting Costas polynomial.
    Without loss of generality, assume $f$ is reduced modulo $x^q-x$, so it has degree at most $q-1$.
    Let 
    \[
        f(x) = \sum_{n=0}^{q-1} c_n x^n.
    \]
    Fix $d\in\F_q$, $d \neq 1$, and take $a\in\F_q^*$ such that
    \[
        f(dx)-f(x) = f(ax).
    \]
    Then, $f(x) = f(dx)-f(ax)$ and by comparing coefficients we must have $c_n = c_n(d^n-a^n)$, for $n=0,1,\dots, q-1$.
    This implies 
    \begin{equation}\label{eq4.1}
        c_n=0 \quad\textup{or}\quad d^n-a^n = 1.
    \end{equation} 
    If for some $n\in\set{0,1,\dots, q-1}$, $c_n \neq 0$, we claim that $\gcd(n,q-1) = 1$.
    
    Let us prove the claim.
    For the sake of a contradiction, assume that for some $n\in\set{0,1,\dots, q-1}$, we have $c_n \neq 0$ and $\gcd(n,q-1)>1$.
    Let $\alpha$ be a primitive element in $\F_q$, and set $d = \alpha^r$, where 
    \[
        r = \frac{q-1}{\gcd(n,q-1)}.
    \]
    Then $rn$ is a mutiple of $q-1$ implying $\alpha^{rn} = 1$.
    By \eqref{eq4.1}, $c_n \neq 0$ implies $d^n-a^n = 1$.
    Hence
    \[
        a^n = d^n-1 = (\alpha^r)^n - 1 = 0,
    \]
    implying $a=0$.
    This is a contradiction because $a\in\F_q^*$.
    We conclude that all the terms of $f$ have $x$ to the power of a number relatively prime to $q-1$.
    
    First, assume $f$ has a linear term.
    Then, by \eqref{eq4.1}, $a = d-1$.
    If $f$ has another term $a_nx^n$ for some $n>1$, $a^n = (d-1)^n$ and \eqref{eq4.1} implies 
    \begin{equation}\label{eq4.2}
        (d-1)^n = d^n-1.
    \end{equation}
    Notice that, if $\F_q$ has odd characteristic, $n$ must be odd because $\gcd(n,q-1) = 1$ and $q-1$ is even.
    Hence, in $\F_q$, $(-1)^n = -1$.
    Otherwise, if $\F_q$ has even characteristic, $(-1)^n = -1 = 1$ for any integer $n$.
    In either case, using \eqref{eq4.2}, $(d-1)^n = d^n + (-1)^n = d^n - 1$.
    Therefore, using the binomial theorem, 
    \begin{align*}
        d^n-1 = \sum_{k=0}^n \binom{n}{k} (-1)^k d^{n-k} =& \; d^n + \sum_{k=1}^{n-1} \binom{n}{k}(-1)^k d^{n-k} + (-1)^n\\
        \implies& \sum_{k=1}^{n-1} \binom{n}{k}(-1)^k d^{n-k} = 0 \quad \forall d\in\F_q, d\neq1.
    \end{align*}
    That is, the polynomial $\sum_{k=1}^{n-1} \binom{n}{k}(-1)^k x^{n-k} \in \F_q[x]$ of degree $n-1\leqslant q-2$ has $q-1$ zeros, so it must be the zero polynomial. 
    Since $\F_q$ has characteristic $p$, we must have $\binom{n}{k}\equiv 0 \pmod{p}$ for $k = 1, 2, \dots n-1$.
    By \cite{fine1947binomial}*{Theorem 3}, $n$ has to be a power of $p$.
    We conclude that $f$ is a linearized permutation polynomial.
    
    Now assume $f$ does not have linear term.
    If $c_s x^s$ is a term of $f$, we showed that $s$ is a unit modulo $q-1$.
    Therefore, there is some $t$ relatively prime to $q-1$, for which $st \equiv 1 \pmod{q-1}$, implying that $f(x^t)$ has a linear term when reduced modulo $x^q-x$.
    Let $L(x)$ be the reduction of $f(x^t)$ modulo $x^q-x$.
    By \Cref{lemma:composition}, $L(x)$ is a shifting Costas polynomial, and, as shown above, it is a linearized permutation polynomial.
    Then, $f(x) \equiv L(x^s) \pmod{x^q-x}$.
    
    Conversely, assume $f(x) \equiv L(x^s)\pmod{x^q-x}$ for some linearized permutation polynomial $L(x)\in\F_q[x]$ and for some $s$, $\gcd(s,q-1) = 1$.
    Notice that $f(x) \equiv L(x^s)\pmod{x^q-x}$ implies that $f$ and $L(x^s)$ are equal as functions over $\F_q$.
    Hence, it is enough to show that $f(x) = L(x^s)$ is a shifting Costas polynomial.
    
    Assume $f(x) = L(x^s)$.
    Clearly $f$ is a permutation polynomial and $f(0) = 0$.
    Fix $d\in\F_q, d\neq 1$.
    Then
    \begin{align*}
        f(dx)-f(x) &= L(d^sx^s)-L(x^s)\\
            &= L((d^s-1)x)
    \end{align*}
    Since $d^s-1\in\F_q$ and $d^s-1\neq 0$, there exists $a\in\F_q^*$ for which $a^s = d^s-1$.
    Therefore, 
    \[
        f(dx)-f(x) = L(a^sx^s) = f(ax).
    \]
    We conclude that $f$ is a shifting Costas polynomial, and the proof is finished.
\end{proof}

By following the proof of \Cref{thm:md_shifting} for polynomials over prime fields, we get a completely different proof for \Cref{prop:shifting} than that by Moreno \cite{moreno1992shifting}*{p.~159}.

\section{Conclusions}

We presented definitions, theory and conjectures that extend results and relations established for (two-dimensional) Costas arrays to  multiple dimensions. Some previous definitions can be seen now as special cases of this theory. 

With this general setting we connected Muratovi\'c-Ribi\'c et al.'s conjecture on Costas polynomials being the composition of linearized permutation polynomials with permutation monomials and a multidimensional version of the Golomb-Moreno conjecture that states that circular Costas maps are multidimensional Welch.
Following Moreno's approach of defining a stronger property to characterize the two-dimensional Welch construction, we defined shifting Costas polynomials to characterize multidimensional Welch maps.
The multidimensional version of the Golomb-Moreno conjecture, that standard circular Costas maps characterize multidimensional Welch maps, remains open.

Also, the connection of circular Costas maps with permutation monomials composed with linearized permutation polynomials allowed us to prove conjectures by Ortiz-Ubarri et al. on bounds on the number of certain circular Costas maps.
We also proved their conjecture on circular Costas maps over elementary abelian groups.

\section{Acknowledgments}
We are very grateful to David Thomson for the conversations and his insight on the generalization of the Golomb-Moreno conjecture which motivated this work.
We also thank Jos\'e Ortiz-Ubarri and Rafael Arce-Nazario for sharing more details on their work in \cite{ortiz2013algebraic} and for their comments on the results presented in \cite{thesis}. This research was funded by the ``Fondo Institucional Para la Investigaci\'on (FIPI)'' from the University of Puerto Rico, R\'{\i}o Piedras.


\bibliographystyle{plain}
\bibliography{ref}

\end{document}